\documentclass[submission%
]{dmtcs}
\usepackage[latin1]{inputenc}
\usepackage[round]{natbib}
\usepackage{amsmath,tikz}
\usepackage{tabulary}
\usepackage{longtable}
\usepackage{tabu}
\usepackage{stmaryrd}
\usepackage{multirow}
\usepackage{multicol}

\let\set\mathbb
\def\sgn{\operatorname{sgn}}
\def\clap#1{\hbox to0pt{\hss#1\hss}}

\def\Stepset#1#2#3#4#5#6#7#8#9{%
  \raisebox{-9pt}{%
  \setlength{\unitlength}{.9pt}%
  \begin{picture}(20,20)(-10,-10)
    \put(-5,-5){\clap{$\scriptstyle\ifx1#1\bullet\else\cdot\fi$}}
    \put(0,-5){\clap{$\scriptstyle\ifx1#2\bullet\else\cdot\fi$}}
    \put(5,-5){\clap{$\scriptstyle\ifx1#3\bullet\else\cdot\fi$}}
    \put(-5,0){\clap{$\scriptstyle\ifx1#4\bullet\else\cdot\fi$}}
    \put(0,0){\clap{$\scriptstyle\ifx1#5\bullet\else\cdot\fi$}}
    \put(5,0){\clap{$\scriptstyle\ifx1#6\bullet\else\cdot\fi$}}
    \put(-5,5){\clap{$\scriptstyle\ifx1#7\bullet\else\cdot\fi$}}
    \put(0,5){\clap{$\scriptstyle\ifx1#8\bullet\else\cdot\fi$}}
    \put(5,5){\clap{$\scriptstyle\ifx1#9\bullet\else\cdot\fi$}}
  \end{picture}}\StepsetB}

\def\StepsetB#1#2#3#4#5#6#7#8{%
  \raisebox{-8pt}{%
  \setlength{\unitlength}{.9pt}%
  \begin{picture}(20,20)(-10,-10)
    \put(-5,-5){\clap{$\scriptstyle\ifx1#1\bullet\else\cdot\fi$}}
    \put(0,-5){\clap{$\scriptstyle\ifx1#2\bullet\else\cdot\fi$}}
    \put(5,-5){\clap{$\scriptstyle\ifx1#3\bullet\else\cdot\fi$}}
    \put(-5,0){\clap{$\scriptstyle\ifx1#4\bullet\else\cdot\fi$}}
    \put(5,0){\clap{$\scriptstyle\ifx1#5\bullet\else\cdot\fi$}}
    \put(-5,5){\clap{$\scriptstyle\ifx1#6\bullet\else\cdot\fi$}}
    \put(0,5){\clap{$\scriptstyle\ifx1#7\bullet\else\cdot\fi$}}
    \put(5,5){\clap{$\scriptstyle\ifx1#8\bullet\else\cdot\fi$}}
  \end{picture}}\StepsetC}

\def\StepsetC#1#2#3#4#5#6#7#8#9{%
  \raisebox{-7pt}{%
  \setlength{\unitlength}{.9pt}%
  \begin{picture}(20,20)(-10,-10)
    \put(-5,-5){\clap{$\scriptstyle\ifx1#1\bullet\else\cdot\fi$}}
    \put(0,-5){\clap{$\scriptstyle\ifx1#2\bullet\else\cdot\fi$}}
    \put(5,-5){\clap{$\scriptstyle\ifx1#3\bullet\else\cdot\fi$}}
    \put(-5,0){\clap{$\scriptstyle\ifx1#4\bullet\else\cdot\fi$}}
    \put(0,0){\clap{$\scriptstyle\ifx1#5\bullet\else\cdot\fi$}}
    \put(5,0){\clap{$\scriptstyle\ifx1#6\bullet\else\cdot\fi$}}
    \put(-5,5){\clap{$\scriptstyle\ifx1#7\bullet\else\cdot\fi$}}
    \put(0,5){\clap{$\scriptstyle\ifx1#8\bullet\else\cdot\fi$}}
    \put(5,5){\clap{$\scriptstyle\ifx1#9\bullet\else\cdot\fi$}}
  \end{picture}}}

\author{Axel Bacher\addressmark{1}\thanks{Email: \email{bacher@lipn.fr}},
  Manuel Kauers\addressmark{2}\thanks{Email: \email{manuel.kauers@jku.at}. Partially supported by the Austria FWF grants Y464-N18 and F50-04.}
  \and Rika Yatchak\addressmark{2}\thanks{Email: \email{rika.yatchak@jku.at}. Partially supported by the Austria FWF grant F50-04}}
\title{Continued Classification of 3D Lattice Models in the Positive~Octant}
\address{\addressmark{1}Laboratoire d'Informatique de Paris Nord, Universit\'{e} Paris 13 \\ 
	 \addressmark{2}Institute for Algebra, Johannes Kepler University, Linz, Austria}
\received{tba}
\revised{tba}
\accepted{tba}
\begin{document}
\maketitle

\begin{abstract}
\paragraph{Abstract.} 
We continue the investigations of lattice walks in the three dimensional lattice restricted to
the positive octant. We separate models which clearly have a D-finite generating function from
models for which there is no reason to expect that their generating function is D-finite, and
we isolate a small set of models whose nature remains unclear and requires further investigation.
For these, we give some
experimental results about their asymptotic behaviour, based on the inspection of a large number
of initial terms. At least for some of them, the guessed asymptotic form seems to tip the balance
towards non-D-finiteness.

\end{abstract}

\section{Introduction}

The past years have seen many contributions to the theory of lattice walks in
the first quadrant~$\set N^2$.  For a fixed set $S\subseteq\set Z^2$, the
quantity of interest is the number $q_{n,i,j}$ of lattice walks of length~$n$
from $(0,0)$ to $(i,j)$ not stepping out of~$\set N^2$. More formally, these
walks can be viewed as elements $(s_1,s_2,\dots,s_n)\in S^n$ with the property
that $\sum_{\ell=1}^m s_\ell\in\set N^2$ for $m=0,\dots,n$ and $\sum_{\ell=1}^n
s_\ell=(i,j)$.  One of the key questions in this context is whether or not the
generating function
$Q(x,y,t)=\sum_{n=0}^\infty\sum_{i=0}^\infty\sum_{j=0}^\infty
q_{n,i,j}x^iy^jt^n$ is D-finite with respect to~$t$.  The answer to this
question depends on~$S$.

For all the sets $S\subseteq\set\{-1,0,1\}^2\setminus\{(0,0)\}$, the nature of
the generating function is known. In their seminal
paper, \cite{bousquet10} reduced the $2^{3^2-1}=256$ different models to 79 among which
there are no obvious bijections and which are not algebraic for classical
reasons~\citep{banderier02,flajolet09}. They then showed that 22 of these 79 models are D-finite
and gave evidence that 56 were not D-finite. The last model required some more
work, but we now have several independent proofs that its generating function is
algebraic, and hence also D-finite~\citep{bostan10,bostan13,bousquet15}. For 51 of the remaining
models, \cite{bostan14a} found that the asymptotic behaviour of the sequence
$q_{n,0,0}$ of walks returning to the origin implies that the power series
$Q(0,0,t)$ is not D-finite (and hence $Q(x,y,t)$ cannot be D-finite
either). Finally, the generating functions of the remaining 5 models are not
D-finite because they have too many singularities~\citep{rechnitzer09,melczer13}.

The case of two dimensions is, in short, well understood. As it turns out, the
D-finiteness of a model is governed by a certain group associated to the model.
There is a general theory~\citep{fayolle99,raschel12,kurkova15} that uniformly explains that the
generating function is D-finite if and only if the associated group is finite.
Among the remaining open problems in 2D are proofs of the non-D-finiteness of
the generating functions $Q(1,1,t)$ for the 56 models where $Q(x,y,t)$ is
non-D-finite, and the classification of models with longer steps, continuing
first results obtained many years ago by \cite{bousquet03}.

What happens in three dimensions? Here, for every fixed
$S\subseteq\set Z^3$ we count the number of walks in $\set N^3$ from $(0,0,0)$
to $(i,j,k)$ with exactly $n$ steps. Again, the principal question is
whether the generating function $Q(x,y,z,t)=\sum_{n=0}^\infty\sum_{i,j,k\in\set
  N} q_{i,j,k,n}x^iy^jz^kt^n$ for the number $q_{i,j,k,n}$ of walks in such a model is
D-finite, and again, the answer depends on the choice
of~$S$. Restricting $S$ to subsets of $\{-1,0,1\}^3\setminus\{(0,0,0)\}$, there
are now $2^{3^3-1}=67{,}108{,}864$ models to be considered.

The classification of these models was initiated by~\cite{bostan09}. More
recently, \cite{bostan14} considered the $\sum_{k=0}^6\binom{26}k=313{,}912$
step sets $S$ with $|S|\leq 6$. After discarding symmetric and simple cases,
they were left with 20,804 models, of which 151 were recognized as D-finite,
and 20,634 were conjectured non-D-finite. No conjecture was made for the remaining
19 cases. A few months ago, \cite{du15} provided non-D-finiteness proofs
for most of the 20,634 models that had been conjectured non-D-finite. The 19
models about which nothing is known have a finite group associated to them,
which, by analogy to the situation in 2D, would suggest that they are
D-finite. On the other hand, \cite{bostan14} were not able to discover any
recurrence or differential equations by guessing, which means that either these
models are not D-finite, or the equations that they satisfy are very large.

In the present paper we continue the classification work for octant walks. We
apply the techniques of \cite{bostan14} to all the 67 million models, and
separate those for which D-finiteness can be proved, those which have no reason
to be D-finite, and those whose nature is unclear. We also give an overview of
the finite groups which appear in~3D; the collection turns out to be a bit more
diverse than in~2D. In the end we found 170 models that are of a
similar nature as the 19 models discovered by \cite{bostan14}. For these we made
an effort to compute a large number of terms and come up with reasonable guesses
for their asymptotic behaviour.

It seems that the exponent $\alpha$ in the asymptotic growth $c\phi^nn^\alpha$
is rational for some cases and irrational for others. This is interesting
because in view of the work of \cite{bostan14a}, an
irrational $\alpha$ can imply non-D-finiteness of the generating functions for
these cases. For those where $\alpha$ seems rational, we have computed
additional terms modulo a prime and used them to try to guess a recurrence or a
differential equation. Several years of computation time has been invested, but
no equation was found. 

\section{Unused Steps, Symmetries, Projectibility, and Decomposibility}

A priori there are $2^{3^3-1}=67{,}108{,}864$ different models. To narrow the number
of cases down to a more manageable size, we apply the same techniques as~\cite{bostan14}.

The \textbf{first filter} sorts out cases that are for simple reasons in bijection to others.
The bijections in question are the permutations of the coordinates. In this filter we
also take care of models containing directions that can never be used. For example,
it is clear that the two models
\Stepset000000000 10100010 000000000
and 
\Stepset000100000 10100010 000000000
consist of exactly the same walks, because no walk of the second model can possibly
involve the step $(-1,0,-1)$ without leaving the first octant.
Here and below, we depict step sets in the same way as~\cite{bostan14}: the first block
contains the directions $(\ast,\ast,-1)$, the second block the directions $(\ast,\ast,0)$,
and the third block the directions $(\ast,\ast,1)$. The diagram on the right therefore
refers to the step set
\[
\{(-1,0,-1),(-1,-1,0),(0,1,0),(1,-1,0)\}.
\]
The number of models surviving this first filter was already worked out by \cite{bostan14} for
step sets with arbitrary cardinalty. According to their Proposition~5, the generating function
$\sum_S u^{|S|}$ where $S$ runs through the essentially different models is given by
\begin{alignat*}1
  &73u^{3}+979u^{4}+6425u^{5}+28071u^{6}+91372u^{7}+234716u^{8}+492168u^{9}\\
  &{}+860382u^{10}+1271488u^{11}+1603184u^{12}+1734396u^{13}+1614372u^{14}\\
  &{}+1293402u^{15}+890395u^{16}+524638u^{17}+263008u^{18}+111251u^{19}\\
  &{}+39256u^{20}+11390u^{21}+2676u^{22}+500u^{23}+73u^{24}+9u^{25}+u^{26}.
\end{alignat*}
Going beyond this counting result, we have written a program that actually lists these models.
The statistics in our listing agrees with this generating function. 

The \textbf{second filter} sorts out models that are in bijection to some lattice
walk model in the quarter plane. For example, for the left model depicted above, all walks
stay within the plane corresponding to the first two coordinates, because there is no step
that moves towards the third dimension. As pointed out by \cite{bostan14}, there are models which
have both steps of type $(\ast,\ast,-1)$ as well as steps of type $(\ast,\ast,1)$, but where
the third dimension nevertheless is immaterial, because the constraint for staying always
nonnegative in the third coordinate is implied by the requirements that the first two
coordinates are nonnegative. Linear programming can be used to detect for a given model
whether it is in bijection to a two-dimensional model. We have applied this test to all the
11,074,225 models passing the first filter and obtained 10,908,263 models that are inherently
three dimensional. The generating function is
\begin{alignat*}1
&u^{3} + 220 u^{4} + 2852 u^{5} + 17731 u^{6} + 70590 u^{7} + 203965 u^{8} + 457650 u^{9}\\
&{}+ 830571 u^{10} + 1251613 u^{11} + 1593013 u^{12} + 1730461 u^{13} + 1613252 u^{14}\\
&{}+ 1293178 u^{15} + 890366 u^{16} + 524636 u^{17} + 263008 u^{18} + 111251 u^{19}\\
&{}+ 39256 u^{20} + 11390 u^{21} + 2676 u^{22} + 500 u^{23} + 73 u^{24}+ 9 u^{25} + u^{26}.
\end{alignat*}
It is fair to eliminate the models that are in bijection with models in the
plane because we have a classification of the generating functions for the
latter. In particular, we have identified~\citep{kauers15c} a list of families
of 2D models with colored steps whose generating function is D-finite. Evidence
was given that this list is complete, at least for models with up to three
colors. These are the models to which a 3D model can possibly be in bijection,
see the article of \cite{bostan14} for more details. 

As a \textbf{third filter,} we determined the models that can be viewed as a direct product
of a one dimensional model and a two dimensional model. These models were called Hadamard
walks by \cite{bostan14}. An example is the decomposition
\[
\Stepset101000101 01000010 101000101
= \overset{\uparrow}{\underset{\downarrow}{\cdot}}
\ \cup \ \left(
\leftarrow{\!\cdot\!}\rightarrow
\ \times \
  \raisebox{-8pt}{%
  \setlength{\unitlength}{.9pt}%
  \begin{picture}(20,20)(-10,-10)
    \put(-5,-5){\clap{$\scriptstyle\bullet$}}
    \put(0,-5){\clap{$\scriptstyle\cdot$}}
    \put(5,-5){\clap{$\scriptstyle\bullet$}}
    \put(-5,0){\clap{$\scriptstyle\cdot$}}
    \put(5,0){\clap{$\scriptstyle\cdot$}}
    \put(-5,5){\clap{$\scriptstyle\bullet$}}
    \put(0,5){\clap{$\scriptstyle\cdot$}}
    \put(5,5){\clap{$\scriptstyle\bullet$}}
  \end{picture}} 
\right).
\]
The generating function for the model on the left can be obtained from the two
generating functions for the models on the right, and since these are D-finite,
the model on the left must also have a D-finite generating function. Conversely,
if a model admits a decomposition into two lower dimensional
models, one of which does not have a D-finite generating function, then this does
not seem to say much about the nature of the generating function for the three
dimensional model. Nevertheless, we decided to discard all the models which
admit a Hadamard decomposition from consideration. We are then left with
alltogether 10,847,434 models, the statistics according to size of step sets
being
\begin{alignat*}1
  &u^{3} + 193u^{4} + 2680u^{5}  + 17238u^{6} + 69542u^{7} + 202072u^{8} + 454485u^{9}\\
  &{} + 826005u^{10}+ 1245615u^{11} + 1585989u^{12} + 1722891u^{13} + 1605940u^{14}\\
  &{}  + 1286692u^{15}  + 885048u^{16}+ 520725u^{17} + 260374u^{18} + 109625u^{19}\\
  &{}+ 38377u^{20} + 10960u^{21} + 2488u^{22} + 436u^{23} + 54u^{24} + 4u^{25}.
\end{alignat*}

\section{The Associated Group}

In their paper on walks in the quarter plane, \cite{bousquet10} make use of a
certain group associated to each model, introduced to the combinatorics
community by~\cite{fayolle99}. \cite{bostan14} consider the following natural analog of
the group for three dimensional models.

For a fixed step set $S\subseteq\{-1,0,1\}^3\setminus\{(0,0,0)\}$,
consider the (Laurent) polynomial
\begin{alignat*}1
  P_S(x,y,z) &= \!\sum_{(u,v,w)\in S} x^u y^v z^w 
  = x^{-1}\kern-.5em\sum_{(-1,v,w)\in S} y^v z^w
  + x^0\kern-.5em\sum_{(0,v,w)\in S} y^v z^w
  + x^1\kern-.5em\sum_{(1,v,w)\in S} y^v z^w.
\end{alignat*}
It is easy to see that $P_S$ remains fixed under the rational transformation
\[
\phi_x\colon\set Q(x,y,z)\to\set Q(x,y,z),\qquad
\phi_x\Bigl(x,y,z\Bigr):=\Bigl(x^{-1} \frac{\sum_{(-1,v,w)\in S} y^v z^w}{\sum_{(1,v,w)\in S} y^v z^w},y,z\Bigr).
\]
In the same way, we can define rational transformations $\phi_y$ and $\phi_z$
which act on $y$ and~$z$, respectively, and leave the other two variables
fixed. The group $G$ associated to the model with step set $S$ is the group
generated by $\phi_x,\phi_y,\phi_z$ under composition.

The generators $\phi_x,\phi_y,\phi_z$ are self-inverse, so 
$G$ can be viewed as a group of words over the alphabet
$\{\phi_x,\phi_y,\phi_z\}$ subject to the relations $\phi_x^2=1$, $\phi_y^2=1$,
$\phi_z^2=1$ and possibly others.  In the case of two dimensions, where the
group only has two generators $\phi_x,\phi_y$, the group is finite if and only
if $(\phi_x\phi_y)^n=1$ for some $n\in\set N$, so the only
finite groups that can appear are the dihedral groups $D_{2n}$ with $2n$
elements.

In the 2D case, to check for the finiteness of such a group it suffices to check whether
some power of $\phi_x\phi_y$ is the identity. In 3D, there is more diversity.
For $n=2,3,\dots$ we consider all the words over $\phi_x,\phi_y,\phi_z$ that are
not equivalent to some shorter word modulo a known relation. If there is no such
word, the group is finite and we stop. Otherwise, for each word in the list, we
check whether the corresponding rational map is the identity. If so, we have found
a new relation and add it to our collection.

As an example, consider the model \Stepset 100000000 00001010 000010000. We have
\[
\phi_x\Bigl(x,y,z\Bigr) = \Bigl(\frac1{xyz},y,z\Bigr),\qquad
\phi_y\Bigl(x,y,z\Bigr) = \Bigl(x,\frac1{xyz},z\Bigr),\qquad
\phi_z\Bigl(x,y,z\Bigr) = \Bigl(x,y,\frac1{xyz}\Bigr).
\]
The only relations of length two are $\phi_x^2=\phi_y^2=\phi_z^2=1$.
There are no new relations of lengths 3, 4, or 5, but we find the relations
$(\phi_x\phi_y)^3=(\phi_x\phi_z)^3=(\phi_y\phi_y)^3=1$ of length~6. These
relations however do not suffice to imply finiteness of the group. Only after we
also find the relations
\begin{alignat*}1
  \phi_z\phi_y\phi_z\phi_x\phi_z\phi_y\phi_z\phi_x =
  \phi_z\phi_y\phi_x\phi_z\phi_x\phi_y\phi_z\phi_x =
  \phi_z\phi_y\phi_z\phi_x\phi_y\phi_z\phi_y\phi_x =
  \phi_y\phi_x\phi_y\phi_z\phi_x\phi_y\phi_x\phi_z = 1
\end{alignat*}
of length~8, it turns out that no word of length 9 can be formed which does not
contain at least one of the found relators as subword. Thus the group is
finite. Its order is~24, the number of words that do not contain any of the
relators as subword.


{}From a set of relations that completely characterizes a finite group, we can
recognize the group using the SmallGroups package in GAP~\citep{GAP4}.
Only 243 of the 10,847,434 models surviving the filters described in the
previous section have a group with $\leq 400$ elements.
This is the \textbf{fourth filter.} The generating function is
\[
8u^{4} + 15u^{6} + 12u^{7} + 21u^{8} + 12u^{9} + 50u^{10} + 24u^{11} + 15u^{12} + 36u^{13} + 20u^{14} + 6u^{15} + 18u^{16} + 6u^{17}.
\]
The groups turn out to be $D_{12}$, $S_4$, $\set Z_2\times S_4$. We believe that
the groups with $>400$ elements are in fact infinite and expect that the
corresponding generating functions are not D-finite. We have also determined the
groups of the 226,791 inherently three dimensional models that admit a Hadamard
decomposition. We found that 2,187 of them have a group with $\leq400$ elements,
these groups are $\set Z_2^3$, $D_{12}$, and $\set Z_2\times D_8$. See
Table~\ref{tab:groups} for more precise statistics. It is noteworthy that we only
encounter Coxeter groups. For example, in the example given above, if we replace
the generator $\phi_z$ by $\psi:=\phi_z\phi_y\phi_z$, then the group is determined
by the relations $\phi_x^2=\phi_y^2=\psi^2=(\phi_x\phi_y)^3=(\phi_x\psi)^2=(\phi_y\psi)^3=1$.

For models with a finite group, we can form the orbit sum $\sum_{g\in
  G} \sgn(g) g(xyz)\in\set Q(x,y,z)$. As explained by \cite{bousquet10} and
\cite{bostan14}, a non-zero orbit sum (almost always) implies D-finiteness of the
generating function $Q(x,y,z)$ of the model. As a \textbf{fifth filter,} we discard
those models which have a finite group but a non-zero orbit sum. This leaves 170
models, the generating function marking step set being
\[
7u^4 + 12u^6 + 8u^7 + 16u^8 + 8u^9 + 35u^{10} + 16u^{11} + 10u^{12} + 24u^{13} + 14u^{14} + 4u^{15} + 12u^{16} + 4u^{17}.
\]
The $7+12=19$ models of with 4 or 6 steps were already identified by~\cite{bostan14}.


\begin{table}
  \begin{center}
  \begin{tabulary}{\textwidth}{|C|C||C|C|}
  \hline
    Group & Hadamard & Non-Hadamard \newline Nonzero O.S. & Non-Hadamard \newline Zero O.S. \\ \hline
    $\mathbb{Z}_{2}\times \mathbb{Z}_{2} \times \mathbb{Z}_{2}$ & 1852 & 0 & 0\\ \hline
    $D_{12}$ & 253 & 66 & 132 \\ \hline
    $\mathbb{Z}_2 \times D_{8}$ & 82 & 0 & 0\\ \hline
    $S_4$ & 0 & 5 & 26\\ \hline
    $\mathbb{Z}_2\times S_4$ & 0 & 2 & 12\\
    \hline 
  \end{tabulary}
  \end{center}
  \caption{Number of models with finite group}
  \label{tab:groups}
\end{table}

\section{Indecomposible and unprojectible models with finite group and zero orbit sum}


We now have a closer look at the 170 non-equivalent models which are not projectible to quarter-plane models,
which do not admit a Hadamard decomposition, whose associated group is finite, and whose orbit sum is zero.
For none of these models, it is known whether the associated generating function is D-finite.
Our goal is to get some idea about the possibilities to be D-finite by looking at the asymptotic behaviour of
the counting sequences for excursions and for walks with arbitrary endpoint, for each of these models. 

\subsection{Computation of Terms}

To compute terms of the generating functions, we used the straightforward
algorithm: we maintain a 3-dimensional array containing all terms
$q_{i,j,k,n}$ for a given value of~$n$. Given the stepset, we can compute from
that the numbers $q_{i,j,k,n+1}$ and iterate over the desired values of~$n$.
Unfortunately, this is extremely costly in both time and memory. We used the
improvements described below to make it more tractable. The C language code for
computing the coefficients was automatically generated from the stepset by a
Sage~\cite{sage} script.

\paragraph{Reduction modulo a prime.} Instead of computing the whole
coefficients, we compute only the residues modulo a prime~$p$.
We chose primes satisfying $p \le 2^{15}$, so that the residues fit in
a 16-bit integer.

\paragraph{Eliminate terms known to be zero.} The next step is to identify the
tuples $(i,j,k,n)$ for which the term~$q_{i,j,k,n}$ can be nonzero and
compute only these terms. From the stepset, we can deduce several inequalities
that have to be satisfied. These inequalities define a 4-dimensional polytope
that we can compute in Sage. Computing only the terms that are in this
polytope saves both time and space.

If we are only interested in the number of excursions of length up to~$N$, we
can further reduce the number of terms that we need: we only need the terms
$q_{i,j,k,n}$ such that there exists an excursion of length~$\le N$ reaching
the point $(i,j,k)$ after~$n$ steps. This adds more constraints to the
polytope.

Finally, in many stepsets, the quadruple $(i,j,k,n)$ has to satisfy modular
constraints (for instance, $i + j + k + n$ has to be even) for $q_{i,j,k,n}$
to be nonzero. This can also be determined automatically. In this case, again,
we do not store the coefficients known to be zero and do not compute them.

\paragraph{Vectorization.} The next optimization is to use the processor's
SIMD instructions (Single Instruction, Multiple Data), that operate on 128-bit
vector registers. Such a register can store eight 16-bit integers in a packed
fashion, and the processor can operate on all of them in parallel with a
single instruction. Moreover, we were able to compute the residues modulo~$p$
without using costly integer divisions. To do that, we note that if
$p\le2^{15}$ and $a$ and $b$ are residues modulo~$p$, their sum modulo~$p$
can be computed as
\def\rem{\operatorname{rem}}
$\rem(a + b, p) = \min(a + b, \rem(a + b - p, 2^{16}))$.
Modular sums can therefore be computed using the vector addition, subtraction,
and minimum instructions present in the SSE4.1 instruction set.

\paragraph{Parallelization.} The final optimization is to distribute the
computations over multiple processors to save time. This is the easiest step:
since all values $q_{i,j,k,n+1}$ can be computed from the values~$q_{i,j,k,n}$
independently of each other, we can give each processor a share of them.
This was done automatically using the OpenMP interface.

\subsection{Guessed Asymptotic Behavior} 

For all the 170 models in question, we calculated the first 2001 terms of the
series $Q(0,0,0,t)$ counting excursions and of the series $Q(1,1,1,t)$ counting
walks with arbitrary endpoint. We want to get an idea about the asymptotic
behaviour of these sequences as $n\to\infty$. We assume that the asymptotic
behaviour of these sequences is given by a linear combination of
terms~$\phi^nn^\alpha$, and our goal is to determine accurate estimates for the
constants $\phi$ and~$\alpha$. For estimating these constants, it is not enough
to know the terms of the counting sequences modulo a prime as computed by the
code described above. However, we can apply the code for several primes and
reconstruct the integer values from the various modular images using Chinese
remaindering. The number of primes needed depends on the size of the integer to
be reconstructed. As an upper bound, we can use that in a model with stepset~$S$
there can be at most $|S|^n$ walks of length~$n$, so if we use primes in the
range $2^{14}\ldots2^{15}$, then $\lceil \frac{n\log|S|}{14\log(2)}\rceil$
primes will always be enough. For our largest step sets, which have 17 elements, we
used 584 primes to recover the 2000th term.

If a sequence $(a_n)_{n=0}^\infty$ grows like $c\phi^nn^\alpha$ for some constants
$\phi,\alpha,c$, then we have $\lim_{n\to\infty} \frac{a_{n+1}}{a_n}=\phi$ and
$\lim_{n\to\infty} \frac{n(a_{n+1}-\phi a_n)}{\phi a_n}=\alpha$. We can thus get
first estimates for $\phi$ and $\alpha$ by simply evaluating the respective
expressions for some large index~$n$. For example, for the counting sequence
$(a_n)_{n=0}^\infty$ of walks with arbitrary endpoint in the model 
\Stepset11011100011100110000110111 
we have
\[ 
  \phi\approx\frac{a_{1200}}{a_{1199}}
  =\frac{25407\text{(\dots 1378 digits\dots)}93695}{17572\text{(\dots 1377 digits\dots)}52363}
  \approx 14.4585690074019.
\]  
Comparison to
\[
  \phi\approx\frac{a_{1190}}{a_{1189}}
  =\frac{63641\text{(\dots 1366 digits\dots)}06567}{44016\text{(\dots 1365 digits\dots)}32175}
  \approx 14.4583480279347
\]
suggests that the accuracy is close to~$10^{-4}$.
To get a better estimate, we use a convergence acceleration technique
due to~\cite{richardson27}. It is based on the assumption that the convergent
sequence $(u_n)_{n=0}^\infty$ whose limit~$u$ is to be estimated has an asymptotic
expansion of the form
\[
  u_n = u + c_1 n^{-1} + c_2 n^{-2} + O(n^{-3})\quad(n\to\infty)
\]
where $c_1,c_2,\dots$ are some unknown constants. By cancellation of the term $c_1n^{-1}$,
it obtains an auxiliary sequence which converges to the same limit but with speed
$1/n^2$ rather than~$1/n$. The cancellation can be achieved in several ways. In
particular, we have
\[
  2 u_{2n} - u_n = u + 0 n^{-1} - \frac12 c_2 n^{-2} + O(n^{-3})\quad(n\to\infty)
\]
and
\[
  (n+1) u_{n+1} - n u_n = u + 0 n^{-1} + c_2 \underbrace{\frac1{n(n+1)}}_{=n^{-2}+O(n^{-3})} + O(n^{-3})\quad(n\to\infty).
\]
Clearly, both versions can be generalized such as to eliminate further terms in the
expansion by taking suitable linear combinations of $u_n,u_{2n},u_{4n},\dots,u_{2^in}$
or $u_n,u_{n+1},\dots,u_{n+i}$, respectively.

For the sequence $(a_n)_{n=0}^\infty$ counting walks with arbitrary endpoints
in the model
\Stepset11011100011100110000110111 
we set $u_n=\frac{a_n}{a_{n-1}}$ and have, for example, 
\begin{alignat*}3
  \frac1{7!}\sum_{k=0}^7(-1)^{k+7}\binom 7k (1200-k)^7 u_{1200-k} &\approx 14.48528121823356265,\\
  \frac1{7!}\sum_{k=0}^7(-1)^{k+7}\binom 7k (1190-k)^7 u_{1190-k} &\approx 14.48528121635317802.
\end{alignat*}
The expressions on the left are such that they cancel the first six terms in the asymptotic
expansion of~$u_n$.
These approximations are accurate enough to recover from them, using standard techniques like
LLL~\citep{lenstra82} or PSLQ~\citep{ferguson92}, that $\phi$ is probably $2(1+\sqrt6)$. 

Once a reliable guess for $\phi$ is available, $\alpha$ can be estimated in the
same way. The results of our calculations are given in tables in the
appendix. It is noteworthy that for some models, $\alpha$ is easily recognized
as a rational number, while for other models, our estimates seem to suggest
that these $\alpha$'s are irrational, which in analogy with~\cite{bostan14a} would imply that the corresponding
sequences are not D-finite. We restrict the focus to the models where the counting sequence for excursions as
well as the counting sequence for walks with arbitrary endpoint have a
(conjecturally) rational~$\alpha$. Of each of these sequences we have
computed the first 5127 terms modulo 16381. (It turned out that 5127 was the
largest number of terms our C code can compute on a machine with 512G of RAM; we
used some twenty such machines, each equipped with 32 processors running in parallel for
several days.)  Regrettably, for none of the sequences we were able to obtain
plausible candidates for potential recurrence or differential equations. If such
equations exist, they must have high order or degree. We are not entirely
convinced that no such equations exist, in view of the example
\Stepset00010000110001000000100001 for which \cite{bostan14} showed that it has a
recurrence of order 55 and degree 3815, whose construction needs some 20000
terms. 


\bibliographystyle{abbrvnat}
\bibliography{sample}

\appendix 
\section{Tables}

\newcount\modelidx
\modelidx=0

\def\one{1}

\newenvironment{mytable}{%
  \par\medskip\noindent
  \def\lhs{1}%
  \hbox to 5mm{\hss idx\hss}%
  \hbox to 2cm{\hss step set\hss}%
  \hbox to 1cm{\hss $(x,y)$\hss}%
  \hbox to 39mm{\hss Asymptotics\hss}%
  \smash{\rule[-1em]{\fboxrule}{2em}}%
  \hbox to 5mm{\hss idx\hss}%
  \hbox to 2cm{\hss step set\hss}%
  \hbox to 1cm{\hss $(x,y)$\hss}%
  \hbox to 39mm{\hss Asymptotics\hss}%
  \hfill\\[-1ex]\rule{\hsize}{\fboxrule}%
}{%
  \par\medskip
}
\def\row#1#2#3#4{%
  \global\advance\modelidx by1\relax
  \vbox to7.74mm{\vfill\hbox to 5mm{\hss\strut\the\modelidx\hss}\vfill}%
  \vbox to7.74mm{\vfill\hbox to 2cm{\hss#2\hss}\vfill}%
  \vbox to7.74mm{\vfill
    \hbox to 1cm{\hss$(0,0)\vphantom{\bigr)^{3n}}$\hss}%
    \hbox to 1cm{\hss$(1,1)\vphantom{\bigr)^{3n}}$\hss}\vfill}%
  \vbox to7.74mm{\vfill
    \hbox to 39mm{\hss$#3\vphantom{\bigr)^{3n}}$\hss}%
    \hbox to 39mm{\hss$#4\vphantom{\bigr)^{3n}}$\hss}\vfill}%
  \ifx\lhs\one
    \def\lhs{0}\smash{\rule{\fboxrule}{12mm}}%
  \else
    \def\lhs{1}\\[-1em]\rule{\hsize}{\fboxrule}
  \fi
}

\scriptsize

\def\left{}
\def\right{}

\begin{mytable}
  \row{o8a08414}{\Stepset00101000001000010000010100}{[n]_{2}\,6^{n} n^{-8.0256624}}{\,6^{n} n^{-3.2634617}}%
\row{o8a03014}{\Stepset00101000000011000000010100}{[n]_{2}\,6^{n} n^{-5.9706049}}{\,6^{n} n^{-2.2353017}\, (*)}%
\row{o9204812}{\Stepset01001000000100100000010010}{[n]_{2}\,6^{n} n^{-5.5631102}}{\,6^{n} n^{-2.0315321}\, (*)}%
\row{o9103022}{\Stepset01000100000011000000100010}{[n]_{2}\,6^{n} n^{-4.5566911}}{\,6^{n} n^{-1.5283424}\, (*)}%
\row{o940840a}{\Stepset01010000001000010000001010}{[n]_{2}\,6^{n} n^{-3.5478909}}{\,6^{n} n^{-1.0239354}\, (*)}%
\row{oa208411}{\Stepset10001000001000010000010001}{[n]_{2}\,6^{n} n^{-3.1240844}}{\,6^{n} n^{-0.8120415}\, (*)}%
\row{o0907824}{\Stepset00100100000111100000100100}{[n]_{2}\,8^{n} n^{-8.0256639}}{\,8^{n} n^{-3.2628309}\, (*)}%
\row{o1b00036}{\Stepset01101100000000000000110110}{[n]_{2}\,8^{n} n^{-5.9706049}}{\,8^{n} n^{-2.2353022}\, (*)}%
\row{o1320132}{\Stepset01001100100000000100110010}{[n]_{2}\,8^{n} n^{-5.5631088}}{\,8^{n} n^{-2.0315369}\, (*)}%
\row{o1e0001e}{\Stepset01111000000000000000011110}{[n]_{2}\,8^{n} n^{-4.5566911}}{\,8^{n} n^{-1.5283262}\, (*)}%
\row{o3300033}{\Stepset11001100000000000000110011}{[n]_{2}\,8^{n} n^{-3.5478909}}{\,8^{n} n^{-1.0239455}\, (*)}%
\row{o360001b}{\Stepset11011000000000000000011011}{[n]_{2}\,8^{n} n^{-3.1240844}}{\,8^{n} n^{-0.8120411}\, (*)}%

\end{mytable}
\begin{center}\footnotesize
  \textbf{\refstepcounter{table}Tab.~\thetable:}
  Models with group $G=\set Z_2\times S_4$.
  The notation $[n]_p$ is meant to be $1$ if $p\mid n$ and $0$ otherwise. 
  For sequences marked with $(\ast)$ the growth seems to be of the form $c(n)\phi^nn^\alpha$
  where $c(n)$ depends on the parity of~$n$. In other cases, the growth looks like $c\phi^nn^\alpha$
  for some constant~$c$.
\end{center}

\medskip
\begin{mytable}
  \row{o0209004}{\Stepset00100000000010010000010000}{[n]_{4}\,4^{n} n^{-5.6432110}}{\,4^{n} n^{-2.0750861}}%
\row{o0802410}{\Stepset00001000001001000000000100}{[n]_{4}\,4^{n} n^{-5.6432110}}{\,4^{n} n^{-2.0591708}}%
\row{o0206002}{\Stepset01000000000001100000010000}{[n]_{4}\,4^{n} n^{-4.7409029}}{\,4^{n} n^{-1.62223}}%
\row{o0404410}{\Stepset00001000001000100000001000}{[n]_{4}\,4^{n} n^{-4.7409029}}{\,4^{n} n^{-1.6186}}%
\row{o0409002}{\Stepset01000000000010010000001000}{[n]_{4}\,4^{n} n^{-3.6508693}}{\,4^{n} n^{-1.075}}%
\row{o020a001}{\Stepset10000000000001010000010000}{[n]_{4}\,4^{n} n^{-3.3257569}}{\,4^{n} n^{-0.9171490}}%
\row{o2001410}{\Stepset00001000001010000000000001}{[n]_{4}\,4^{n} n^{-3.3257569}}{\,4^{n} n^{-0.91}}%
\row{o9104822}{\Stepset01000100000100100000100010}{[n]_{2}\,6^{n} n^{-5.6432112}}{\,6^{n} n^{-2.0716051}\, (*)}%
\row{o9083042}{\Stepset01000010000011000001000010}{[n]_{2}\,6^{n} n^{-4.7409028}}{\,6^{n} n^{-1.6204494}\, (*)}%
\row{oa204811}{\Stepset10001000000100100000010001}{[n]_{2}\,6^{n} n^{-3.6508694}}{\,6^{n} n^{-1.0754}}%
\row{o941020a}{\Stepset01010000010000001000001010}{[n]_{2}\,6^{n} n^{-3.3257568}}{\,6^{n} n^{-0.9128785}\, (*)}%
\row{o0a0b414}{\Stepset00101000001011010000010100}{[n]_{2}\,8^{n} n^{-5.6432107}}{\,8^{n} n^{-2.0716056}\, (*)}%
\row{o1207812}{\Stepset01001000000111100000010010}{[n]_{2}\,8^{n} n^{-4.7409030}}{\,8^{n} n^{-1.6204500}\, (*)}%
\row{o140b40a}{\Stepset01010000001011010000001010}{[n]_{2}\,8^{n} n^{-3.6508692}}{\,8^{n} n^{-1.0754330}\, (*)}%
\row{o220b411}{\Stepset10001000001011010000010001}{[n]_{2}\,8^{n} n^{-3.3257570}}{\,8^{n} n^{-0.9128784}\, (*)}%
\row{o9906c32}{\Stepset01001100001101100000100110}{\,10^{n} n^{-5.6432111}}{\,10^{n} n^{-2.0716055}}%
\row{o930d826}{\Stepset01100100000110110000110010}{\,10^{n} n^{-5.6432111}}{\,10^{n} n^{-2.0716054}}%
\row{o9487452}{\Stepset01001010001011100001001010}{\,10^{n} n^{-4.7409028}}{\,10^{n} n^{-1.6204517}}%
\row{o8e0e40e}{\Stepset01110000001001110000011100}{\,10^{n} n^{-4.7409028}}{\,10^{n} n^{-1.6204516}}%
\row{oa60d813}{\Stepset11001000000110110000011001}{\,10^{n} n^{-3.6508693}}{\,10^{n} n^{-1.0754347}}%
\row{o961a20b}{\Stepset11010000010001011000011010}{\,10^{n} n^{-3.3257569}}{\,10^{n} n^{-0.9128784}}%
\row{ob41161a}{\Stepset01011000011010001000001011}{\,10^{n} n^{-3.3257569}}{\,10^{n} n^{-0.9128784}}%
\row{o9b0fc36}{\Stepset01101100001111110000110110}{\,14^{n} n^{-5.6432110}}{\,14^{n} n^{-2.0716054}}%
\row{o9e0fc1e}{\Stepset01111000001111110000011110}{\,14^{n} n^{-4.7409029}}{\,14^{n} n^{-1.6204519}}%
\row{ob60fc1b}{\Stepset11011000001111110000011011}{\,14^{n} n^{-3.6508693}}{\,14^{n} n^{-1.0754348}}%
\row{ob61b61b}{\Stepset11011000011011011000011011}{\,14^{n} n^{-3.3257569}}{\,14^{n} n^{-0.9128784}}%

\end{mytable}
\begin{center}\footnotesize
  \textbf{\refstepcounter{table}Tab.~\thetable:}
  Models with group $G=S_4$. Same notational conventions as in the previous table.
\end{center}

\medskip
\begin{mytable}
  \row{oa1400a1}{\Stepset10000101000000000010100001}{[n]_{2}\,6^{n} n^{-4}}{\,6^{n} n^{-5/4}\, (*)}%
\row{o942010a}{\Stepset01010000100000000100001010}{[n]_{2}\,6^{n} n^{-4}}{\,6^{n} n^{-5/4}\, (*)}%
\row{oc20a30a}{\Stepset01010000110001010000010000}{[n]_{3}\,(2 \, \sqrt{3} + 3)^n n^{-4}}{\,(2 \, \sqrt{3} + 3)^n n^{-2.25}}%
\row{od106411}{\Stepset10001000001001100000100010}{[n]_{3}\,(2 \, \sqrt{3} + 3)^n n^{-4}}{\,(2 \, \sqrt{3} + 3)^n n^{-2.24}}%
\row{oe1400b1}{\Stepset10001101000000000010100001}{[n]_{2}\,(4 \, \sqrt{3})^n n^{-4}}{\,(4 \, \sqrt{3})^n n^{-11/4}\, (*)}%
\row{od42011a}{\Stepset01011000100000000100001010}{[n]_{2}\,(4 \, \sqrt{3})^n n^{-4}}{\,(4 \, \sqrt{3})^n n^{-11/4}\, (*)}%
\row{oe3400a1}{\Stepset10000101000000000010110001}{[n]_{2}\,(4 \, \sqrt{3})^n n^{-4}}{\,7^{n} n^{-3/4}\, (*)}%
\row{od62010a}{\Stepset01010000100000000100011010}{[n]_{2}\,(4 \, \sqrt{3})^n n^{-4}}{\,7^{n} n^{-3/4}\, (*)}%
\row{oe209822}{\Stepset01000100000110010000010001}{[n]_{3}\,(2 \, \sqrt{3} + 3)^n n^{-4}}{\,7^{n} n^{-3/4}}%
\row{od411414}{\Stepset00101000001010001000001010}{[n]_{3}\,(2 \, \sqrt{3} + 3)^n n^{-4}}{\,7^{n} n^{-3/4}}%
\row{o02114b1}{\Stepset10001101001010001000010000}{\,7^{n} n^{-4}}{\,7^{n} n^{-9/4}}%
\row{o020a31a}{\Stepset01011000110001010000010000}{\,7^{n} n^{-4}}{\,7^{n} n^{-9/4}}%
\row{o23400b1}{\Stepset10001101000000000010110001}{[n]_{2}\,8^{n} n^{-4}}{\,8^{n} n^{-5/4}\, (*)}%
\row{o162011a}{\Stepset01011000100000000100011010}{[n]_{2}\,8^{n} n^{-4}}{\,8^{n} n^{-5/4}\, (*)}%
\row{o220b822}{\Stepset01000100000111010000010001}{\,7^{n} n^{-4}}{\,8^{n} n^{-3/4}}%
\row{o1413414}{\Stepset00101000001011001000001010}{\,7^{n} n^{-4}}{\,8^{n} n^{-3/4}}%
\row{o54201ab}{\Stepset11010101100000000100001010}{[n]_{2}\,(6 \, \sqrt{2})^n n^{-4}}{\,(6 \, \sqrt{2})^n n^{-11/4}\, (*)}%
\row{o61401ab}{\Stepset11010101100000000010100001}{[n]_{2}\,(6 \, \sqrt{2})^n n^{-4}}{\,(6 \, \sqrt{2})^n n^{-11/4}\, (*)}%
\row{o630c823}{\Stepset11000100000100110000110001}{\,9^{n} n^{-4}}{\,9^{n} n^{-5/4}}%
\row{o7104c31}{\Stepset10001100001100100000100011}{\,9^{n} n^{-4}}{\,9^{n} n^{-5/4}}%
\row{o4e1820e}{\Stepset01110000010000011000011100}{\,9^{n} n^{-4}}{\,9^{n} n^{-5/4}}%
\row{o543150a}{\Stepset01010000101010001100001010}{\,9^{n} n^{-4}}{\,9^{n} n^{-5/4}}%
\row{o6514825}{\Stepset10100100000100101000101001}{[n]_{2}\,(6 \, \sqrt{2})^n n^{-4}}{\,9^{n} n^{-3/4}}%
\row{o756010a}{\Stepset01010000100000000110101011}{[n]_{2}\,(6 \, \sqrt{2})^n n^{-4}}{\,9^{n} n^{-3/4}}%
\row{o82115ab}{\Stepset11010101101010001000010000}{\,(2 \, \sqrt{6} + 3)^n n^{-4}}{\,(2 \, \sqrt{6} + 3)^n n^{-2.2502770}}%
\row{o820a3ab}{\Stepset11010101110001010000010000}{\,(2 \, \sqrt{6} + 3)^n n^{-4}}{\,(2 \, \sqrt{6} + 3)^n n^{-2.2502612}}%
\row{o94201bb}{\Stepset11011101100000000100001010}{[n]_{2}\,(2 \, \sqrt{21})^n n^{-4}}{\,(2 \, \sqrt{21})^n n^{-11/4}\, (*)}%
\row{oa1401bb}{\Stepset11011101100000000010100001}{[n]_{2}\,(2 \, \sqrt{21})^n n^{-4}}{\,(2 \, \sqrt{21})^n n^{-11/4}\, (*)}%
\row{o9a0a61a}{\Stepset01011000011001010000010110}{\,(2 \, \sqrt{3} + 6)^n n^{-4}}{\,(2 \, \sqrt{3} + 6)^n n^{-2.2530695}}%
\row{o930e413}{\Stepset11001000001001110000110010}{\,(2 \, \sqrt{3} + 6)^n n^{-4}}{\,(2 \, \sqrt{3} + 6)^n n^{-2.2530524}}%
\row{o96201ab}{\Stepset11010101100000000100011010}{[n]_{2}\,(4 \, \sqrt{6})^n n^{-4}}{\,(4 \, \sqrt{6})^n n^{-11/4}\, (*)}%
\row{oa3401ab}{\Stepset11010101100000000010110001}{[n]_{2}\,(4 \, \sqrt{6})^n n^{-4}}{\,(4 \, \sqrt{6})^n n^{-11/4}\, (*)}%
\row{o8e1920e}{\Stepset01110000010010011000011100}{\,(4 \, \sqrt{3} + 3)^n n^{-4}}{\,(4 \, \sqrt{3} + 3)^n n^{-2.3138908}}%
\row{o943151a}{\Stepset01011000101010001100001010}{\,(4 \, \sqrt{3} + 3)^n n^{-4}}{\,(4 \, \sqrt{3} + 3)^n n^{-2.3139180}}%
\row{oa14a31a}{\Stepset01011000110001010010100001}{\,(4 \, \sqrt{3} + 3)^n n^{-4}}{\,(4 \, \sqrt{3} + 3)^n n^{-2.3138703}}%
\row{o9515415}{\Stepset10101000001010101000101010}{\,(4 \, \sqrt{3} + 3)^n n^{-4}}{\,(4 \, \sqrt{3} + 3)^n n^{-2.3139010}}%
\row{oa30d823}{\Stepset11000100000110110000110001}{\,(4 \, \sqrt{3} + 3)^n n^{-4}}{\,(4 \, \sqrt{3} + 3)^n n^{-2.3138550}}%
\row{oa1514b1}{\Stepset10001101001010001010100001}{\,(4 \, \sqrt{3} + 3)^n n^{-4}}{\,(4 \, \sqrt{3} + 3)^n n^{-2.3138788}}%
\row{ob56011a}{\Stepset01011000100000000110101011}{[n]_{2}\,(4 \, \sqrt{6})^n n^{-4}}{\,10^{n} n^{-3/4}}%
\row{oa515825}{\Stepset10100100000110101000101001}{[n]_{2}\,(4 \, \sqrt{6})^n n^{-4}}{\,10^{n} n^{-3/4}}%
\row{oa30e823}{\Stepset11000100000101110000110001}{\,(4 \, \sqrt{3} + 3)^n n^{-4}}{\,10^{n} n^{-3/4}}%
\row{oa619826}{\Stepset01100100000110011000011001}{\,(2 \, \sqrt{6} + 3)^n n^{-4}}{\,10^{n} n^{-3/4}}%
\row{o9619416}{\Stepset01101000001010011000011010}{\,(2 \, \sqrt{3} + 6)^n n^{-4}}{\,10^{n} n^{-3/4}}%
\row{ob209c32}{\Stepset01001100001110010000010011}{\,(2 \, \sqrt{3} + 6)^n n^{-4}}{\,10^{n} n^{-3/4}}%
\row{oa516825}{\Stepset10100100000101101000101001}{[n]_{2}\,(2 \, \sqrt{21})^n n^{-4}}{\,10^{n} n^{-3/4}}%
\row{ob411c34}{\Stepset00101100001110001000001011}{\,(2 \, \sqrt{6} + 3)^n n^{-4}}{\,10^{n} n^{-3/4}}%
\row{ob76010a}{\Stepset01010000100000000110111011}{[n]_{2}\,(2 \, \sqrt{21})^n n^{-4}}{\,10^{n} n^{-3/4}}%
\row{o9516415}{\Stepset10101000001001101000101010}{\,(4 \, \sqrt{3} + 3)^n n^{-4}}{\,10^{n} n^{-0.7501789}}%
\row{o963150a}{\Stepset01010000101010001100011010}{\,(4 \, \sqrt{3} + 3)^n n^{-4}}{\,10^{n} n^{-0.7486131}}%
\row{oaa0aa2a}{\Stepset01010100010101010000010101}{\,(4 \, \sqrt{3} + 3)^n n^{-4}}{\,10^{n} n^{-0.7484231}}%
\row{ob106c31}{\Stepset10001100001101100000100011}{\,(4 \, \sqrt{3} + 3)^n n^{-4}}{\,10^{n} n^{-0.7481541}}%
\row{o8e1a20e}{\Stepset01110000010001011000011100}{\,(4 \, \sqrt{3} + 3)^n n^{-4}}{\,10^{n} n^{-0.7504432}}%
\row{oc2115bb}{\Stepset11011101101010001000010000}{\,(2 \, \sqrt{7} + 3)^n n^{-4}}{\,(2 \, \sqrt{7} + 3)^n n^{-2.2502209}}%
\row{oc20a3bb}{\Stepset11011101110001010000010000}{\,(2 \, \sqrt{7} + 3)^n n^{-4}}{\,(2 \, \sqrt{7} + 3)^n n^{-2.2502015}}%
\row{od30f413}{\Stepset11001000001011110000110010}{\,10^{n} n^{-4}}{\,10^{n} n^{-9/4}}%
\row{oc21b71a}{\Stepset01011000111011011000010000}{\,10^{n} n^{-4}}{\,10^{n} n^{-9/4}}%
\row{od6201bb}{\Stepset11011101100000000100011010}{[n]_{2}\,(4 \, \sqrt{7})^n n^{-4}}{\,(4 \, \sqrt{7})^n n^{-11/4}\, (*)}%
\row{oe3401bb}{\Stepset11011101100000000010110001}{[n]_{2}\,(4 \, \sqrt{7})^n n^{-4}}{\,(4 \, \sqrt{7})^n n^{-11/4}\, (*)}%
\row{oe30f823}{\Stepset11000100000111110000110001}{\,11^{n} n^{-4}}{\,11^{n} n^{-5/4}}%
\row{oe3514b1}{\Stepset10001101001010001010110001}{\,11^{n} n^{-4}}{\,11^{n} n^{-5/4}}%
\row{oce1b20e}{\Stepset01110000010011011000011100}{\,11^{n} n^{-4}}{\,11^{n} n^{-5/4}}%
\row{od63151a}{\Stepset01011000101010001100011010}{\,11^{n} n^{-4}}{\,11^{n} n^{-5/4}}%
\row{oe517825}{\Stepset10100100000111101000101001}{[n]_{2}\,(4 \, \sqrt{7})^n n^{-4}}{\,11^{n} n^{-3/4}}%
\row{od61b416}{\Stepset01101000001011011000011010}{\,10^{n} n^{-4}}{\,11^{n} n^{-3/4}}%
\row{oe61b826}{\Stepset01100100000111011000011001}{\,(2 \, \sqrt{7} + 3)^n n^{-4}}{\,11^{n} n^{-3/4}}%
\row{of76011a}{\Stepset01011000100000000110111011}{[n]_{2}\,(4 \, \sqrt{7})^n n^{-4}}{\,11^{n} n^{-3/4}}%
\row{of20bc32}{\Stepset01001100001111010000010011}{\,10^{n} n^{-4}}{\,11^{n} n^{-3/4}}%
\row{of413c34}{\Stepset00101100001111001000001011}{\,(2 \, \sqrt{7} + 3)^n n^{-4}}{\,11^{n} n^{-3/4}}%
\row{o14315ab}{\Stepset11010101101010001100001010}{\,(6 \, \sqrt{2} + 3)^n n^{-4}}{\,(6 \, \sqrt{2} + 3)^n n^{-2.2546223}}%
\row{o0f1c20f}{\Stepset11110000010000111000111100}{\,(6 \, \sqrt{2} + 3)^n n^{-4}}{\,(6 \, \sqrt{2} + 3)^n n^{-2.2546131}}%
\row{o21515ab}{\Stepset11010101101010001010100001}{\,(6 \, \sqrt{2} + 3)^n n^{-4}}{\,(6 \, \sqrt{2} + 3)^n n^{-2.2546063}}%
\row{o214a3ab}{\Stepset11010101110001010010100001}{\,(6 \, \sqrt{2} + 3)^n n^{-4}}{\,(6 \, \sqrt{2} + 3)^n n^{-2.2545976}}%
\row{o330cc33}{\Stepset11001100001100110000110011}{\,12^{n} n^{-4}}{\,12^{n} n^{-5/4}}%
\row{o1e1861e}{\Stepset01111000011000011000011110}{\,12^{n} n^{-4}}{\,12^{n} n^{-5/4}}%
\row{o271c827}{\Stepset11100100000100111000111001}{\,(6 \, \sqrt{2} + 3)^n n^{-4}}{\,12^{n} n^{-3/4}}%
\row{o3514c35}{\Stepset10101100001100101000101011}{\,(6 \, \sqrt{2} + 3)^n n^{-4}}{\,12^{n} n^{-3/4}}%
\row{o2e18a2e}{\Stepset01110100010100011000011101}{\,(6 \, \sqrt{2} + 3)^n n^{-4}}{\,12^{n} n^{-3/4}}%
\row{o357150a}{\Stepset01010000101010001110101011}{\,(6 \, \sqrt{2} + 3)^n n^{-4}}{\,12^{n} n^{-3/4}}%
\row{o54315bb}{\Stepset11011101101010001100001010}{\,(2 \, \sqrt{21} + 3)^n n^{-4}}{\,(2 \, \sqrt{21} + 3)^n n^{-2.2523243}}%
\row{o61515bb}{\Stepset11011101101010001010100001}{\,(2 \, \sqrt{21} + 3)^n n^{-4}}{\,(2 \, \sqrt{21} + 3)^n n^{-2.2523199}}%
\row{o4f1d20f}{\Stepset11110000010010111000111100}{\,(2 \, \sqrt{21} + 3)^n n^{-4}}{\,(2 \, \sqrt{21} + 3)^n n^{-2.2523167}}%
\row{o614a3bb}{\Stepset11011101110001010010100001}{\,(2 \, \sqrt{21} + 3)^n n^{-4}}{\,(2 \, \sqrt{21} + 3)^n n^{-2.2523126}}%
\row{o56315ab}{\Stepset11010101101010001100011010}{\,(4 \, \sqrt{6} + 3)^n n^{-4}}{\,(4 \, \sqrt{6} + 3)^n n^{-2.2721995}}%
\row{o4f1e20f}{\Stepset11110000010001111000111100}{\,(4 \, \sqrt{6} + 3)^n n^{-4}}{\,(4 \, \sqrt{6} + 3)^n n^{-2.2721868}}%
\row{o63515ab}{\Stepset11010101101010001010110001}{\,(4 \, \sqrt{6} + 3)^n n^{-4}}{\,(4 \, \sqrt{6} + 3)^n n^{-2.2721586}}%
\row{o634a3ab}{\Stepset11010101110001010010110001}{\,(4 \, \sqrt{6} + 3)^n n^{-4}}{\,(4 \, \sqrt{6} + 3)^n n^{-2.2721471}}%
\row{o5e1961e}{\Stepset01111000011010011000011110}{\,(4 \, \sqrt{3} + 6)^n n^{-4}}{\,(4 \, \sqrt{3} + 6)^n n^{-2.3215931}}%
\row{o571d417}{\Stepset11101000001010111000111010}{\,(4 \, \sqrt{3} + 6)^n n^{-4}}{\,(4 \, \sqrt{3} + 6)^n n^{-2.3215822}}%
\row{o7a09e3a}{\Stepset01011100011110010000010111}{\,(4 \, \sqrt{3} + 6)^n n^{-4}}{\,(4 \, \sqrt{3} + 6)^n n^{-2.3215670}}%
\row{o730dc33}{\Stepset11001100001110110000110011}{\,(4 \, \sqrt{3} + 6)^n n^{-4}}{\,(4 \, \sqrt{3} + 6)^n n^{-2.3215560}}%
\row{o571e417}{\Stepset11101000001001111000111010}{\,(4 \, \sqrt{3} + 6)^n n^{-4}}{\,13^{n} n^{-3/4}}%
\row{o671e827}{\Stepset11100100000101111000111001}{\,(2 \, \sqrt{21} + 3)^n n^{-4}}{\,13^{n} n^{-3/4}}%
\row{o7515c35}{\Stepset10101100001110101000101011}{\,(4 \, \sqrt{6} + 3)^n n^{-4}}{\,13^{n} n^{-3/4}}%
\row{o7516c35}{\Stepset10101100001101101000101011}{\,(2 \, \sqrt{21} + 3)^n n^{-4}}{\,13^{n} n^{-3/4}}%
\row{o6e1aa2e}{\Stepset01110100010101011000011101}{\,(2 \, \sqrt{21} + 3)^n n^{-4}}{\,13^{n} n^{-3/4}}%
\row{o6e19a2e}{\Stepset01110100010110011000011101}{\,(4 \, \sqrt{6} + 3)^n n^{-4}}{\,13^{n} n^{-3/4}}%
\row{o757151a}{\Stepset01011000101010001110101011}{\,(4 \, \sqrt{6} + 3)^n n^{-4}}{\,13^{n} n^{-3/4}}%
\row{o777150a}{\Stepset01010000101010001110111011}{\,(2 \, \sqrt{21} + 3)^n n^{-4}}{\,13^{n} n^{-3/4}}%
\row{o5e1a61e}{\Stepset01111000011001011000011110}{\,(4 \, \sqrt{3} + 6)^n n^{-4}}{\,13^{n} n^{-0.6802255}}%
\row{o730ec33}{\Stepset11001100001101110000110011}{\,(4 \, \sqrt{3} + 6)^n n^{-4}}{\,13^{n} n^{-0.6795702}}%
\row{o7a0ae3a}{\Stepset01011100011101010000010111}{\,(4 \, \sqrt{3} + 6)^n n^{-4}}{\,13^{n} n^{-0.6785484}}%
\row{o671d827}{\Stepset11100100000110111000111001}{\,(4 \, \sqrt{6} + 3)^n n^{-4}}{\,13^{n} n^{-0.7499959}}%
\row{o96315bb}{\Stepset11011101101010001100011010}{\,(4 \, \sqrt{7} + 3)^n n^{-4}}{\,(4 \, \sqrt{7} + 3)^n n^{-2.2583992}}%
\row{oa3515bb}{\Stepset11011101101010001010110001}{\,(4 \, \sqrt{7} + 3)^n n^{-4}}{\,(4 \, \sqrt{7} + 3)^n n^{-2.2583717}}%
\row{oa34a3bb}{\Stepset11011101110001010010110001}{\,(4 \, \sqrt{7} + 3)^n n^{-4}}{\,(4 \, \sqrt{7} + 3)^n n^{-2.2583622}}%
\row{o8f1f20f}{\Stepset11110000010011111000111100}{\,(4 \, \sqrt{7} + 3)^n n^{-4}}{\,(4 \, \sqrt{7} + 3)^n n^{-2.2583888}}%
\row{ob30fc33}{\Stepset11001100001111110000110011}{\,14^{n} n^{-4}}{\,14^{n} n^{-5/4}}%
\row{o9e1b61e}{\Stepset01111000011011011000011110}{\,14^{n} n^{-4}}{\,14^{n} n^{-5/4}}%
\row{oa71f827}{\Stepset11100100000111111000111001}{\,(4 \, \sqrt{7} + 3)^n n^{-4}}{\,14^{n} n^{-3/4}}%
\row{ob517c35}{\Stepset10101100001111101000101011}{\,(4 \, \sqrt{7} + 3)^n n^{-4}}{\,14^{n} n^{-3/4}}%
\row{oae1ba2e}{\Stepset01110100010111011000011101}{\,(4 \, \sqrt{7} + 3)^n n^{-4}}{\,14^{n} n^{-3/4}}%
\row{ob77151a}{\Stepset01011000101010001110111011}{\,(4 \, \sqrt{7} + 3)^n n^{-4}}{\,14^{n} n^{-3/4}}%
\row{odf1c61f}{\Stepset11111000011000111000111110}{\,(6 \, \sqrt{2} + 6)^n n^{-4}}{\,(6 \, \sqrt{2} + 6)^n n^{-2.2567151}}%
\row{ofb0ce3b}{\Stepset11011100011100110000110111}{\,(6 \, \sqrt{2} + 6)^n n^{-4}}{\,(6 \, \sqrt{2} + 6)^n n^{-2.2566922}}%
\row{of71cc37}{\Stepset11101100001100111000111011}{\,(6 \, \sqrt{2} + 6)^n n^{-4}}{\,15^{n} n^{-3/4}}%
\row{ofe18e3e}{\Stepset01111100011100011000011111}{\,(6 \, \sqrt{2} + 6)^n n^{-4}}{\,15^{n} n^{-3/4}}%
\row{o1f1d61f}{\Stepset11111000011010111000111110}{\,(2 \, \sqrt{21} + 6)^n n^{-4}}{\,(2 \, \sqrt{21} + 6)^n n^{-2.2533372}}%
\row{o215b7bb}{\Stepset11011101111011011010100001}{\,(2 \, \sqrt{21} + 6)^n n^{-4}}{\,(2 \, \sqrt{21} + 6)^n n^{-2.2533261}}%
\row{o3b0ee3b}{\Stepset11011100011101110000110111}{\,(4 \, \sqrt{6} + 6)^n n^{-4}}{\,(2 \, \sqrt{42} + 3)^n n^{-2.2792897}}%
\row{o1f1e61f}{\Stepset11111000011001111000111110}{\,(4 \, \sqrt{6} + 6)^n n^{-4}}{\,(4 \, \sqrt{6} + 6)^n n^{-2.2793307}}%
\row{o2f1da2f}{\Stepset11110100010110111000111101}{\,(2 \, \sqrt{42} + 3)^n n^{-4}}{\,(2 \, \sqrt{42} + 3)^n n^{-2.3617}}%
\row{o35715bb}{\Stepset11011101101010001110101011}{\,(2 \, \sqrt{42} + 3)^n n^{-4}}{\,(2 \, \sqrt{42} + 3)^n n^{-2.15}}%
\row{o371dc37}{\Stepset11101100001110111000111011}{\,(4 \, \sqrt{6} + 6)^n n^{-4}}{\,16^{n} n^{-3/4}}%
\row{o371ec37}{\Stepset11101100001101111000111011}{\,(2 \, \sqrt{21} + 6)^n n^{-4}}{\,16^{n} n^{-3/4}}%
\row{o3e19e3e}{\Stepset01111100011110011000011111}{\,(4 \, \sqrt{6} + 6)^n n^{-4}}{\,16^{n} n^{-3/4}}%
\row{o3e1ae3e}{\Stepset01111100011101011000011111}{\,(2 \, \sqrt{21} + 6)^n n^{-4}}{\,16^{n} n^{-3/4}}%
\row{o2f1ea2f}{\Stepset11110100010101111000111101}{\,(2 \, \sqrt{42} + 3)^n n^{-4}}{\,16^{n} n^{-0.7859779}}%
\row{o37715ab}{\Stepset11010101101010001110111011}{\,(2 \, \sqrt{42} + 3)^n n^{-4}}{\,16^{n} n^{-0.7859706}}%
\row{o5f1f61f}{\Stepset11111000011011111000111110}{\,(4 \, \sqrt{7} + 6)^n n^{-4}}{\,(4 \, \sqrt{7} + 6)^n n^{-2.2614300}}%
\row{o635b7bb}{\Stepset11011101111011011010110001}{\,(4 \, \sqrt{7} + 6)^n n^{-4}}{\,(4 \, \sqrt{7} + 6)^n n^{-2.2613989}}%
\row{o7e1be3e}{\Stepset01111100011111011000011111}{\,(4 \, \sqrt{7} + 6)^n n^{-4}}{\,17^{n} n^{-3/4}}%
\row{o771fc37}{\Stepset11101100001111111000111011}{\,(4 \, \sqrt{7} + 6)^n n^{-4}}{\,17^{n} n^{-3/4}}%

\end{mytable}
\begin{center}\footnotesize
  \textbf{\refstepcounter{table}Tab.~\thetable:}
  Models with group $G=D_{12}$. Same notational conventions as in the previous table.
\end{center}

\end{document}